\documentclass[12pt]{amsart}
\usepackage{amsfonts}
\usepackage{graphicx}
\usepackage{amsmath}
\usepackage{amssymb}
\input{amssym.def}
\newtheorem{theorem}{Theorem}

\newtheorem{lemma}{Lemma}
\newtheorem{remark}{Remark}

\theoremstyle{remark}

\newcommand{\im}{\text{\rm Im }}

\input epsf
\begin{document}
\title[L\"owner equation]{Exponential driving function for the L\"owner equation}
\author[D.~Prokhorov]{Dmitri Prokhorov}

\address{D.~Prokhorov: Department of Mathematics and Mechanics, Saratov State University, Saratov
410012, Russia} \email{ProkhorovDV@info.sgu.ru}

\begin{abstract}
We consider the chordal L\"owner differential equation with the model driving function $\root3\of
t$. Holomorphic and singular solutions are represented by their series. It is shown that a
disposition of values of different singular and branching solutions is monotonic, and solutions to
the L\"owner equation map slit domains onto the upper half-plane. The slit is a $C^1$-curve. We
give an asymptotic estimate for the ratio of harmonic measures of the two slit sides.
\end{abstract}
\maketitle

\section{Introduction}

The L\"owner differential equation introduced by K.~L\"owner \cite{Loewner} served a source to
study properties of univalent functions on the unit disk. Nowadays it is of growing interest in
many areas, see, e.g., \cite{Markina}. The L\"owner equation for the upper half-plane $\mathbb H$
appeared later (see, e.g., \cite{Aleksandrov}) and became popular during the last decades. Define a
function $w=f(z,t)$, $z\in\mathbb H$, $t\geq0$,
\begin{equation}
f(z,t)=z+\frac{2t}{z}+O\left(\frac{1}{z^2}\right), \;\;\; z\to\infty, \label{exp}
\end{equation}
which maps $\mathbb H\setminus K_t$ onto $\mathbb H$ and solves the {\it chordal} L\"owner ordinary
differential equation
\begin{equation}
\frac{df(z,t)}{dt}=\frac{2}{f(z,t)-\lambda(t)},\;\;\;f(z,0)=z,\;\;\;z\in\mathbb H, \label{Leo}
\end{equation}
where the driving function $\lambda(t)$ is continuous and real-valued.

The conformal maps $f(z,t)$ are continuously extended onto $z\in\mathbb R$ minus the closure of
$K_t$ and the extended map also satisfies equation (\ref{Leo}). Following \cite{LMR}, we pay
attention to an old problem to determine, in terms of $\lambda$, when $K_t$ is a Jordan arc,
$K_t=\gamma(t)$, $t\geq0$, emanating from the real axis $\mathbb R$. In this case $f(z,t)$ are
continuously extended onto the two sides of $\gamma(t)$,
\begin{equation}
\lambda(t)=f(\gamma(t),t),\;\;\;\gamma(t)=f^{-1}(\lambda(t),t). \label{gam}
\end{equation}

Points $\gamma(t)$ are treated as prime ends which are different for the two sides of the arc. Note
that Kufarev \cite{Kufarev} proposed a counterexample of the non-slit mapping for the {\it radial}
L\"owner equation in the disk. For the chordal L\"owner equation, Kufarev's example corresponds to
$\lambda(t)=3\sqrt2\sqrt{1-t}$, see \cite{Kager}, \cite{LMR} for details.

Equation (\ref{Leo}) admits integrating in quadratures for partial cases of $\lambda(t)$ studied in
\cite{Kager}, \cite{PrZ}. The integrability cases of (\ref{Leo}) are invariant under linear and
scaling transformations of $\lambda(t)$, see, e.g., \cite{LMR}. Therefore, assume without loss of
generality that $\lambda(0)=0$ and, equivalently, $\gamma(0)=0$.

The picture of singularity lines for driving functions $\lambda(t)$ belonging to the Lipschitz
class $\text{Lip}(1/2)$ with the exponent 1/2 is well studied, see, e.g., \cite{LMR} and references
therein.

This article is aimed to show that in the case of the cubic root driving function
$\lambda(t)=\root3\of t$ in (\ref{Leo}), that is,
\begin{equation}
\frac{df(z,t)}{dt}=\frac{2}{f(z,t)-\root3\of t},\;\;\;f(z,0)=z,\;\;\;\im z\geq0, \label{Le3}
\end{equation}
the solution $w=f(z,t)$ is a slit mapping for $t>0$ small enough, i.e., $K_t=\gamma(t)$, $0<t<T$.

The driving function $\lambda(t)=\root3\of t$ is chosen as a typical function of the Lipschitz
class $\text{Lip}(1/3)$. We do not try to cover the most general case but hope that the model
driving function serves a demonstration for a wider class. By the way, the case when the trace
$\gamma$ is a circular arc meeting the real axis tangentially is studied in \cite{ProkhVas}. The
explicit solution for the inverse function gave a driving term of the form $\lambda(t)=C
t^{1\over3}+\dots$ which corresponds to the above driving function asymptotically.

The main result of the article is contained in the following theorem which shows that $f(z,t)$ is a
mapping from a slit domain $D(t)=\mathbb H\setminus\gamma(t)$.

\begin{theorem}
Let $f(z,t)$ be a solution to the L\"owner equation (\ref{Le3}). Then $f(\cdot,t)$ maps
$D(t)=\mathbb H\setminus\gamma(t)$ onto $\mathbb H$ for $t>0$ small enough where $\gamma(t)$ is a
$C^1$-curve, except probably for the point $\gamma(0)=0$.
\end{theorem}

Preliminary results of Section 2 in the article concern the theory of differential equations and
preparations for the main proof.

Theorem 1 together with helpful lemmas are proved in Section 3.

Section 4 is devoted to estimates for harmonic measures of the two sides of the slit generated by
the L\"owner equation (\ref{Le3}). Theorem 2 in this Section gives the asymptotic relation for the
ratio of these harmonic measures as $t\to0$.

In Section 5 we consider holomorphic solutions to (\ref{Le3}) represented by power series and
propose asymptotic expansions for the radius of convergence of the series.

\section{Preliminary statements}

Change variables $t\to\tau^3$, $g(z,\tau):=f(z,\tau^3)$, and reduce equation (\ref{Le3}) to

\begin{equation}
\frac{dg(z,\tau)}{d\tau}=\frac{6\tau^2}{g(z,\tau)-\tau}, \;\;\;g(z,0)=z, \;\;\;\im z\geq0.
\label{Le2}
\end{equation}

Note that differential equations $$\frac{dy}{dx}=\frac{Q(x,y)}{P(x,y)}$$ with holomorphic functions
$P(x,y)$ and $Q(x,y)$ are well known both for complex and real variables, especially in the case of
polynomials $P$ and $Q$, see, e.g., \cite{Bendixson}, \cite{Borel}, \cite{Golubev},
\cite{Poincare}, \cite{Sansone-1}, \cite{Sansone-2}.

If $z\neq0$, then $g(z,0)\neq0$, and there exists a {\it regular} solution $g(z,\tau)$ to
(\ref{Le2}) holomorphic in $\tau$ for $|\tau|$ small enough which is unique for every $z\neq0$. We
are interested mostly in studying {\it singular} solutions to (\ref{Le2}), i.e., those which do not
satisfy the uniqueness conditions for equation (\ref{Le2}). Every point $(g(z_0,\tau_0),\tau_0)$
such that $g(z_0,\tau_0)=\tau_0$ is a {\it singular point} for equation (\ref{Le2}). If
$\tau_0\neq0$, then $(g(z_0,\tau_0),\tau_0)$ is an {\it algebraic solution critical point}, and
corresponding singular solutions to (\ref{Le2}) through this point are expanded in series in terms
$(\tau-\tau_0)^{1/m}$, $m\in\mathbb N$. So these singular solutions are different branches of the
same analytic function, see [17,~Chap.9, \S1].

The point $(g(z_0,\tau_0),\tau_0)=(0,0)$ is the only {\it singular point of indefinite character}
for (\ref{Le2}). It is determined when the numerator and denominator in the right-hand side of
(\ref{Le2}) vanish simultaneously. All the singular solutions to (\ref{Le2}) which are not branches
of the same analytic function pass through this point $(0,0)$ [17,~Chap.9, \S1].

Regular and singular solutions to (\ref{Le2}) behave according to the Poincar\'e-Bendixson theorems
\cite{Poincare}, \cite{Bendixson}, [17,~Chap.9, \S1]. Namely, two integral curves of differential
equation (\ref{Le2}) intersect only at the singular point $(0,0)$. An integral curve of (\ref{Le2})
can have multiple points only at $(0,0)$. Bendixson \cite{Bendixson} considered real integral
curves globally and stated that they have endpoints at knots and focuses and have an extension
through a saddle. Under these assumptions, the Bendixson theorem \cite{Bendixson} makes possible
only three cases for equation (\ref{Le2}) in a neighborhood of $(0,0)$: (a) an integral curve is
closed, i.e., it is a cycle; (b) an integral curve is a spiral which tends to a cycle
asymptotically; (c) an integral curve has the endpoint at $(0,0)$.

Recall the integrability case \cite{Kager} of the L\"owner differential equation (\ref{Leo}) with
the square root forcing $\lambda(t)=c\sqrt t$. After changing variables $t\to\tau^2$, the singular
point $(0,0)$ in this case is a saddle according to the Poincar\'e classification \cite{Poincare}
for linear differential equations. From the other side, another integrability case \cite{Kager}
with the square root forcing $\lambda(t)=c\sqrt{1-t}$, after changing variables $t\to1-\tau^2$,
leads to the focus at $(0,0)$.

Going back to equation (\ref{Le2}) remark that its solutions are infinitely differentiable with
respect to the real variable $\tau$, see [4,~Chap.1, \S1], [17,~Chap.9, \S1]. Hence recurrent
evaluations of Taylor coefficients can help to find singular solutions provided that a resulting
series will have a positive convergence radius [16,~Chap.3, \S1]. Apply this method to equation
(\ref{Le2}). Let
\begin{equation}
g_s(0,\tau)=\sum_{n=1}^{\infty}a_n\tau^n \label{Si1}
\end{equation}
be a a formal power series for singular solutions to (\ref{Le2}). Note that $g_s$ is not
necessarily unique. It depends on the path along which $z$ approaches to 0, $z\notin K_{\tau}$.
Substitute (\ref{Si1}) into (\ref{Le2}) and see that
\begin{equation}
\sum_{n=1}^{\infty}na_n\tau^{n-1}\left(\sum_{n=1}^{\infty}a_n\tau^n-\tau\right)=6\tau^2.
\label{Si2}
\end{equation}
Equating coefficients at the same powers in both sides of (\ref{Si2}) obtain that
\begin{equation}
a_1(a_1-1)=0. \label{Si3}
\end{equation}

This equation gives two possible values $a_1=1$ and $a_1=0$ to two singular solutions $g^+(0,\tau)$
and $g^-(0,\tau)$. In both cases equation (\ref{Si2}) implies recurrent formulas for coefficients
$a_n^+$ and $a_n^-$ of $g^+(0,\tau)$ and $g^-(0,\tau)$ respectively,
\begin{equation}
a_1^+=1,\;\;a_2^+=6,\;\;a_n^+=-\sum_{k=2}^{n-1}ka_k^+a_{n+1-k}^+, \;\;n\geq3, \label{Si4}
\end{equation}

\begin{equation}
a_1^-=0,\;\;a_2^-=-3,\;\; a_n^-=\frac{1}{n}\sum_{k=2}^{n-1}ka_k^-a_{n+1-k}^-,\;\;n\geq3,
\label{Si5}
\end{equation}

Show that the series $\sum_{n=1}^{\infty}a_n^+\tau^n$ formally representing $g^+(0,\tau)$ diverges
for all $\tau\neq0$.

\begin{lemma}
For $n\geq2$, the inequalities
\begin{equation}
6^{n-1}(n-1)!\leq|a_n^+|\leq12^{n-1}n^{n-3} \label{Si6}
\end{equation}
hold.
\end{lemma}

\proof For $n=2$, the estimate (\ref{Si6}) from below holds with the equality sign. Suppose that
these estimates are true for $k=2,\dots,n-1$ and substitute them in (\ref{Si2}). For $n\geq3$, we
have $$|a_n^+|=\sum_{k=2}^{n-1}k|a_k^+||a_{n+1-k}^+|
\geq\sum_{k=2}^{n-1}k6^{k-1}(k-1)!6^{n-k}(n-k)!=$$
$$6^{n-1}\sum_{k=2}^{n-1}k!(n-k)!\geq6^{n-1}(n-1)!\;.$$ This confirms by induction the estimate
(\ref{Si6}) from below.

Similarly, for $n=2,3$, the estimate (\ref{Si6}) from above is easily verified. Suppose that these
estimates are true for $k=2,\dots,n-1$ and substitute them in (\ref{Si2}). For $n\geq4$, we have
$$|a_n^+|=\sum_{k=2}^{n-1}k|a_k^+||a_{n+1-k}^+|
\leq\sum_{k=2}^{n-1}k12^{k-1}k^{k-3}12^{n-k}(n+1-k)^{n-2-k}=$$
$$12^{n-1}\sum_{k=2}^{n-1}k^{k-2}(n+1-k)^{n-2-k}
\leq12^{n-1}\left(\sum_{k=2}^{n-2}(n-1)^{k-2}(n-1)^{n-2-k}+\frac{(n-1)^{n-3}}{2}\right)$$
$$<12^{n-1}\left(\sum_{k=2}^{n-2}(n-1)^{n-4}+(n-1)^{n-4}\right)<12^{n-1}n^{n-3}$$ which completes
the proof.
\endproof

Evidently, the upper estimates (\ref{Si6}) are preserved for $|a_n^-|$, $n\geq2$.

The lower estimates (\ref{Si6}) imply divergence of $\sum_{n=1}^{\infty}a_n^+\tau^n$ for
$\tau\neq0$. Therefore equation (\ref{Le2}) does not have any holomorphic solution in a
neighborhood of $\tau_0=0$. There exist some methods to summarize the series
$\sum_{n=1}^{\infty}a_n^+\tau^n$, the Borel regular method among them \cite{Borel}, [16,~Chap.3,
\S1]. Let $$G(\tau)=\sum_{n=1}^{\infty}\frac{a_n^+}{n!}\tau^n,$$ this series converges for
$|\tau|<1/12$ according to Lemma 1. The Borel sum equals
$$h(\tau)=\int_0^{\infty}e^{-x}G(\tau x)dx$$ and solves (\ref{Le2}) provided it determines an
analytic function. The same approach is applied to $\sum_{n=1}^{\infty}a_n^-\tau^n$.

In any case solutions $g_1(0,\tau)$, $g_2(0,\tau)$ to (\ref{Le2}) emanating from the singular point
$(0,0)$ satisfy the asymptotic relations $$g_1(0,\tau)=\sum_{k=1}^na_k^+\tau^k+o(\tau^n),\;\;\;
g_2(0,\tau)=\sum_{k=1}^na_k^-\tau^k+o(\tau^n),\;\;\;\tau\to0,$$ for all $n\geq2$, $o(\tau^n)$ in
both representations depend on $n$.

Let $f_1(0,t):=g_1(0,\tau^3)$, $f_2(0,t):=g_2(0,\tau^3)$. Since $f_1(0,t)=\root3\of
t+6\root3\of{t^2}+o(\root3\of{t^2})$ and $f_2(0,t)=-3\root3\of{t^2}+o(\root3\of{t^2})$ as $t\to0$,
the inequality $$f_2(0,t)<\root3\of t<f_1(0,t)$$ holds for all $t>0$ small enough.

Let us find representations for all other singular solutions to equation (\ref{Le3}) which appear
at $t>0$. Suppose there is $z_0\in\mathbb H$ and $t_0>0$ such that $f(z_0,t_0)=\root3\of t$. Then
$(f(z_0,t_0),t_0)$ is a singular point of equation (\ref{Le3}), and $f(z_0,t)$ is expanded in
series with powers $(t-t_0)^{n/m}$, $m\in\mathbb N$,

\begin{equation}
f(z_0,t)=\root3\of t_0+\sum_{n=1}^{\infty}b_{n/m}(t-t_0)^{n/m}. \label{Si7}
\end{equation}

Substitute (\ref{Si7}) into (\ref{Le3}) and see that
$$\sum_{n=1}^{\infty}\frac{nb_{n/m}(t-t_0)^{n/m-1}}{m}\times$$
\begin{equation}
\left(\sum_{n=1}^{\infty}b_{n/m}
(t-t_0)^{n/m}-\sum_{n=1}^{\infty}\frac{(-1)^{n-1}2\cdot5\dots(3n-4)}{n!}
\frac{(t-t_0)^n}{(3t_0)^n}\right)=2. \label{Si8}
\end{equation}

Equating coefficients at the same powers in both sides of (\ref{Si8}) obtain that $m=2$ and
\begin{equation}
(b_{1/2})^2=4. \label{Si9}
\end{equation}
This equation gives two possible values $b_{1/2}=2$ and $b_{1/2}=-2$ to two branches $f_1(z_0,t)$
and $f_2(z_0,t)$ of the solution (\ref{Si7}). Indeed, we can accept only one of possibilities, for
example $b_{1/2}=2$, while the second case is obtained by going to another branch of
$(t-t_0)^{n/2}$ when passing through $t=t_0$. So we have recurrent formulas for coefficients
$b_{n/2}$ of $f_1(z_0,t)$ and $f_2(z_0,t)$,

\begin{equation}
b_{1/2}=2,\;\;b_{n/2}=\frac{1}{n+1}\left(c_{n/2}-\frac{1}{2}\sum_{k=2}^{n-1}kb_{k/2}
(b_{(n+1-k)/2}-c_{(n+1-k)/2})\right),\;\;n\geq2, \label{Si10}
\end{equation}

where
\begin{equation}
c_{(2k-1)/2}=0,\;\;c_k=\frac{(-1)^{k-1}2\cdot5\dots(3k-4)}{3^kt_0^{k-1/3}k!},\;\; k=1,2,\dots\;.
\label{Si11}
\end{equation}

Since
$$f_1(z_0,t)=\root3\of{t_0}+2\sqrt{t-t_0}+o(\sqrt{t-t_0}),\;
f_2(z_0,t)=\root3\of{t_0}-2\sqrt{t-t_0}+o(\sqrt{t-t_0}),$$
$$\root3\of t=
\root3\of{t_0}+\frac{1}{3t_0}(t-t_0)+o(t-t_0),\;\;\;t\to t_0+0,$$ the inequality
$$f_2(z_0,t)<\root3\of t<f_1(z_0,t)$$ holds for all $t>t_0$ close to $t_0$.

\section{Proof of the main results}

The theory of differential equations claims that integral curves of equation (\ref{Le3}) intersect
only at the singular point $(0,0)$ [17,~Chap.9, \S1]. In particular, this implies the local
inequalities $f_2(0,t)<f_2(z_0,t)<\root3\of t<f_1(z_0,t)<f_1(0,t)$ where $(f(z_0,t_0),t_0)$ is an
algebraic solution critical point for equation (\ref{Le3}). We will give an independent proof of
these inequalities which can be useful for more general driving functions.

\begin{lemma}
For $t>0$ small enough and a singular point $(f(z_0,t_0),t_0)$ for equation (\ref{Le3}), $0<t_0<t$,
the following inequalities
$$f_2(0,t)<f_2(z_0,t)<\root3\of t<f_1(z_0,t)<f_1(0,t)$$ hold.
\end{lemma}

\proof To show that $f_1(z_0,t)<f_1(0,t)$ let us subtract equations
$$\frac{df_1(0,t)}{dt}=\frac{2}{f_1(0,t)-\root3\of t},\;\;\;f_1(0,0)=0,$$
$$\frac{df_1(z_0,t)}{dt}=\frac{2}{f_1(z_0,t)-\root3\of t}, \;\;\;f_1(z_0,t_0)=\root3\of{t_0},$$
and obtain

$$\frac{d(f_1(0,t)-f_1(z_0,t))}{dt}=\frac{2(f_1(z_0,t)-f_1(0,t))}
{(f_1(0,t)-\root3\of t)(f_1(z_0,t)-\root3\of t)},$$ which can be written in the form
$$\frac{d\log(f_1(0,t)-f_1(z_0,t))}{dt}=\frac{-2}
{(f_1(0,t)-\root3\of t)(f_1(z_0,t)-\root3\of t)}.$$ Suppose that $T>t_0$ is the smallest number for
which $f_1(0,T)=f_1(z_0,T)$. This implies that

\begin{equation}
\int_{t_0}^T\frac{dt}{(f_1(0,t)-\root3\of t)(f_1(z_0,t)-\root3\of t)}= \infty. \label{Si12}
\end{equation}

To evaluate the integral in (\ref{Si12}) we should study the behavior of $f_1(z_0,t)-\root3\of t$
with the help of differential equation

\begin{equation}
\frac{d(f_1(z_0,t)-\root3\of t)}{dt}= \frac{2}{f_1(z_0,t)-\root3\of
t}-\frac{1}{3\root3\of{t^2}}=\frac{\root3\of t+6\root3\of{t^2}-f_1(z_0,t)}{3\root3\of{t^2}
(f_1(z_0,t)-\root3\of t)}. \label{Si13}
\end{equation}

Calculate that $a_3^+=-72$ and write the asymptotic relation
$$f_1(0,t)=\root3\of t+6\root3\of{t^2}-72t+o(t),\;\;t\to+0.$$ There exists a number $T'>0$ such
that for $0<t<T'$, $\root3\of t+6\root3\of{t^2}>f_1(0,t)$. Consequently, the right-hand side in
(\ref{Si13}) is positive for $0<t<T'$. Note that $T'$ does not depend on $t_0$. The condition
"$t>0$ small enough" in Lemma 2 is understood from now as $0<t<T'$. We see from (\ref{Si13}) that
for such $t$, $f_1(z_0,t)-\root3\of t$ is increasing with $t$, $t_0<t<T<T'$. Therefore, the
integral in the left-hand side of (\ref{Si12}) is finite. The contradiction against equality
(\ref{Si12}) denies the existence of $T$ with the prescribed properties which proves the third and
the fourth inequalities in Lemma 2.

The rest of inequalities in Lemma 2 are proved similarly and even easier. To show that
$f_2(z_0,t)>f_2(0,t)$ let us subtract equations
$$\frac{df_2(0,t)}{dt}=\frac{2}{f_2(0,t)-\root3\of t},\;\;\;f_2(0,0)=0,$$
$$\frac{df_2(z_0,t)}{dt}=\frac{2}{f_2(z_0,t)-\root3\of t},
\;\;\;f_2(z_0,t_0)=\root3\of{t_0},$$ and obtain

$$\frac{d(f_2(0,t)-f_2(z_0,t))}{dt}=\frac{2(f_2(z_0,t)-f_2(0,t))}
{(f_2(0,t)-\root3\of t)(f_2(z_0,t)-\root3\of t)},$$ which can be written in the form
$$\frac{d\log(f_2(z_0,t)-f_2(0,t))}{dt}=\frac{-2}
{(f_2(0,t)-\root3\of t)(f_2(z_0,t)-\root3\of t)}.$$ Suppose that $T>t_0$ is the smallest number for
which $f_2(z_0,T)=f_2(0,T)$. This implies that

\begin{equation}
\int_{t_0}^T\frac{dt}{(f_2(0,t)-\root3\of t)(f_2(z_0,t)-\root3\of t)}= \infty. \label{Si14}
\end{equation}

To evaluate the integral in (\ref{Si14}) we should study the behavior of $f_2(z_0,t)-\root3\of t$
with the help of differential equation

\begin{equation}
\frac{d(f_2(z_0,t)-\root3\of t)}{dt}= \frac{2}{f_2(z_0,t)-\root3\of
t}-\frac{1}{3\root3\of{t^2}}=\frac{\root3\of t+6\root3\of{t^2}-f_2(z_0,t)}{3\root3\of{t^2}
(f_2(z_0,t)-\root3\of t)}. \label{Si15}
\end{equation}

Since
$$f_2(0,t)=-3\root3\of{t^2}+o(\root3\of{t^2}),\;\;t\to+0,$$ there exists a
number $T''>0$ such that for $0<t<T''$, $\root3\of t+6\root3\of{t^2}>f_2(0,t)$. Consequently, the
right-hand side in (\ref{Si15}) is positive for $0<t<T''$. We see from (\ref{Si15}) that for such
$t$, $f_2(0,t)-\root3\of t$ is decreasing with $t$, $t_0<t<T<T''$. Therefore, the integral in the
left-hand side of (\ref{Si14}) is finite. The contradiction against equality (\ref{Si14}) denies
the existence of $T$ with the prescribed properties which completes the proof.
\endproof

Add and complete the inequalities of Lemma 2 by the following statements demonstrating a monotonic
disposition of values for different singular solutions.

\begin{lemma}
For $t>0$ small enough and singular points $(f(z_1,t_1),t_1)$, $(f(z_0,t_0),t_0)$ for equation
(\ref{Le3}), $0<t_1<t_0<t$, the following inequalities
$$f_2(z_1,t)<f_2(z_0,t),\;\;\;f_1(z_0,t)<f_1(z_1,t)$$ hold.
\end{lemma}

\proof Similarly to Lemma 2, subtract equations
$$\frac{df_1(z_1,t)}{dt}=\frac{2}{f_1(z_1,t)-\root3\of t},
\;\;\;f_1(z_1,t_1)=\root3\of{t_1},$$
$$\frac{df_1(z_0,t)}{dt}=\frac{2}{f_1(z_0,t)-\root3\of t},
\;\;\;f_1(z_0,t_0)=\root3\of{t_0},$$ and obtain

$$\frac{d(f_1(z_1,t)-f_1(z_0,t))}{dt}=
\frac{2(f_1(z_0,t)-f_1(z_1,t))}{(f_1(z_1,t)-\root3\of t)(f_1(z_0,t)-\root3\of t)},$$ which can be
written in the form
$$\frac{d\log(f_1(z_1,t)-f_1(z_0,t))}{dt}=\frac{-2}
{(f_1(z_1,t)-\root3\of t)(f_1(z_0,t)-\root3\of t)}.$$ Suppose that $T>t_0$ is the smallest number
for which $f_1(z_1,T)=f_1(z_0,T)$. This implies that

\begin{equation}
\int_{t_0}^T\frac{dt}{(f_1(z_1,t)-\root3\of t)(f_1(z_0,t)-\root3\of t)}= \infty. \label{Si16}
\end{equation}

To evaluate the integral in (\ref{Si16}) apply to (\ref{Si13}) and obtain that there exists a
number $T'>0$ such that for $0<t<T'$, $f_1(z_0,t)-\root3\of t$ is increasing with $t$,
$t_0<t<T<T'$. Therefore, the integral in the left-hand side of (\ref{Si16}) is finite. The
contradiction against equality (\ref{Si16}) denies the existence of $T$ with the prescribed
properties which proves the second inequality of Lemma 3.

To prove the first inequality of Lemma 3 subtract equations
$$\frac{df_2(z_1,t)}{dt}=\frac{2}{f_2(z_1,t)-\root3\of t},\;\;\;
f_2(z_1,t_1)=\root3\of{t_1},$$
$$\frac{df_2(z_0,t)}{dt}=\frac{2}{f_2(z_0,t)-\root3\of t},
\;\;\;f_2(z_0,t_0)=\root3\of{t_0},$$ and obtain after dividing by $f_2(z_1,t)-f_2(z_0,t)$

$$\frac{d\log(f_2(z_0,t)-f_2(z_1,t))}{dt}=\frac{-2}
{(f_2(z_1,t)-\root3\of t)(f_2(z_0,t)-\root3\of t)}.$$ Suppose that $T>t_0$ is the smallest number
for which $f_2(z_0,T)=f_2(z_1,T)$. This implies that

\begin{equation}
\int_{t_0}^T\frac{dt}{(f_2(z_1,t)-\root3\of t)(f_2(z_0,t)-\root3\of t)}=\infty. \label{Si17}
\end{equation}

To evaluate the integral in (\ref{Si17}) apply to (\ref{Si15}) and obtain that $\root3\of
t+6\root3\of{t^2}>\root3\of t>f_2(0,t)$. Consequently, the right-hand side in (\ref{Si15}) is
positive and we see that $f_2(0,t)-\root3\of t$ is decreasing with $t$, $t_0<t<T$. Therefore, the
integral in the left-hand side of (\ref{Si17}) is finite. The contradiction against equality
(\ref{Si17}) denies the existence of $T$ with the prescribed properties which completes the proof.
\endproof

{\it Proof of Theorem 1.}

For $t_0>0$, there is a hull $K_{t_0}\subset\mathbb H$ such that $f(\cdot,t_0)$ maps $\mathbb
H\setminus K_{t_0}$ onto $\mathbb H$. We refer to \cite{LMR} for definitions and more details. The
hull $K_{t_0}$ is driven by $\root3\of t$. The function $f(\cdot,t_0)$ is extended continuously
onto the set of prime ends on $\partial(\mathbb H\setminus K_{t_0})$ and maps this set onto
$\mathbb R$. One of the prime ends is mapped on $\root3\of{t_0}$. Let $z_0=z_0(t_0)$ represent this
prime end.

Lemmas 2 and 3 describe the structure of the pre-image of $\mathbb H$ under $f(\cdot,t)$. All the
singular solutions $f_1(0,t)$, $f_2(0,t)$, $f_1(z_0,t)$, $f_2(z_0,t)$, $0<t_0<t<T'$, are
real-valued and satisfy the inequalities of Lemmas 2 and 3. So the segment $I=[f_2(0,t),f_1(0,t)]$
is the union of the segments $I_2=[f_2(0,t),\root3\of t]$ and $I_1=[\root3\of t,f_1(0,t)]$. The
segment $I_2$ consists of points $f_2(z(\tau),t)$, $0\leq\tau<t$, and the segment $I_1$ consists of
points $f_1(z(\tau),t)$, $0\leq\tau<t$. All these points belong to the boundary $\mathbb
R=\partial\mathbb H$. This means that all the points $z(\tau)$, $0\leq\tau<t$, belong to the
boundary $\partial(\mathbb H\setminus K_t)$ of $\mathbb H\setminus K_t$. Moreover, every point
$z(\tau)$ except for the tip determines exactly two prime ends corresponding to $f_1(z(\tau),t)$
and $f_2(z(\tau),t)$. Evidently, $z(\tau)$ is continuous on $[0,t]$. This proves that
$z(\tau):=\gamma(\tau)$ represents a curve $\gamma:=K_t$ with prime ends corresponding to points on
different sides of $\gamma$. This proves that $f^{-1}(w,t)$ maps $\mathbb H$ onto the slit domain
$\mathbb H\setminus\gamma(t)$ for $t>0$ small enough.

It remains to show that $\gamma(t)$ is a $C^1$-curve. Fix $t_0>0$ from a neighborhood of $t=0$.
Denote $g(w,t)=f^{-1}(w,t)$ an inverse of $f(z,t)$, and $h(w,t):=f(g(w,t_0),t)$, $t\geq t_0$. The
arc $\gamma[t_0,t]:=K_t\setminus K_{t_0}$ is mapped by $f(z,t_0)$ onto a curve $\gamma_1(t)$ in
$\mathbb H$ emanating from $\root3\of{t_0}\in\mathbb R$. So the function $h(w,t)$ is well defined
on $\mathbb H\setminus\gamma_1(t_0)$, $t\geq t_0$. Expand $h(w,t)$ near infinity,
$$h(w,t)=g(w,t_0)+\frac{2t}{g(w,t_0)}+O\left(\frac{1}
{g^2(w,t_0)}\right)=w+\frac{2(t-t_0)}{w}+O\left(\frac{1}{w^2}\right).$$ Such expansion satisfies
(\ref{exp}) after changing variables $t\to t-t_0$. The function $h(w,t)$ satisfies the differential
equation
$$\frac{dh(w,t)}{dt}=\frac{2} {h(w,t)-\root3\of
t},\;\;\;h(w,t_0)=w,\;\;\;w\in\mathbb H.$$ This equation becomes the L\"owner differential equation
if $t_1:=t-t_0$, $h_1(w,t_1):=h(w,t_0+t_1)$,

\begin{equation}
\frac{dh_1(w,t_1)}{dt_1}=\frac{2}{h_1(w,t_1)-\root3\of{t_1+t_0}},\;\;\;
h_1(w,0)=w,\;\;\;w\in\mathbb H. \label{Cu1}
\end{equation}

The driving function $\lambda(t_1)=\root3\of{t_1+t_0}$ in (\ref{Cu1}) is analytic for $t_1\geq0$.
It is known [1,~p.59] that under this condition $h_1(w,t_1)$ maps $\mathbb H\setminus\gamma_1$ onto
$\mathbb H$ where $\gamma_1$ is a $C^1$-curve in $\mathbb H$ emanating from
$\lambda(0)=\root3\of{t_0}$. The same does the function $h(w,t)$.

Go back to $f(z,t)=h(f(z,t_0),t)$ and see that $f(z,t)$ maps $\mathbb H\setminus\gamma(t)$ onto
$\mathbb H$, $\gamma(t)=\gamma[0,t_0]\cup\gamma[t_0,t]$, and $\gamma[t_0,t]$ is a $C^1$-curve.
Tending $t_0$ to 0 we prove that $\gamma(t)$ is a $C^1$-curve, except probably for the point
$\gamma(0)=0$. This completes the proof.

\section{Harmonic measures of the slit sides}

The function $f(z,t)$ solving (\ref{Le3}) maps $\mathbb H\setminus\gamma(t)$ onto $\mathbb H$. The
curve $\gamma(t)$ has two sides. Denote $\gamma_1=\gamma_1(t)$ the side of $\gamma$ which is mapped
by the extended function $f(z,t)$ onto $I_1=[\root3\of t,f_1(0,t)]$. Similarly,
$\gamma_2=\gamma_2(t)$ is the side of $\gamma$ which is the pre-image of $I_2=[f_2(0,t),\root3\of
t]$ under $f(z,t)$.

Remind that the harmonic measures $\omega(f^{-1}(i,t);\gamma_k,\mathbb H\setminus\gamma(t),t)$ of
$\gamma_k$ at $f^{-1}(i,t)$ with respect to $\mathbb H\setminus\gamma(t)$ are defined by the
functions $\omega_k$ which are harmonic on $\mathbb H\setminus\gamma(t)$ and continuously extended
on its closure except for the endpoints of $\gamma$, $\omega_k|_{\gamma_k(t)}=1$,
$\omega_k|_{\mathbb R\cup(\gamma(t)\setminus\gamma_k(t))}=0$, $k=1,2$, see, e.g., [6,~Chap.3,
\S3.6]. Denote
$$m_k(t):=\omega(f^{-1}(i,t);\gamma_k,\mathbb H\setminus\gamma(t),t),\;\;\;k=1,2.$$

\begin{theorem}
Let $f(z,t)$ be a solution to the L\"owner equation (\ref{Le3}). Then
\begin{equation}
\lim_{t\to+0}\frac{m_1(t)}{m_2^2(t)}=6\pi. \label{har}
\end{equation}
\end{theorem}

\proof The harmonic measure is invariant under conformal transformations. So
$$\omega(f^{-1}(i,t);\gamma_k,\mathbb H\setminus\gamma(t),t)=\Omega(i;f(\gamma_k,t),\mathbb H,t)$$
are defined by the harmonic functions $\Omega_k$ which are harmonic on $\mathbb H$ and continuously
extended on $\mathbb R$ except for the endpoints of $f(\gamma_k,t)$, $\Omega_k|_{f(\gamma_k,t)}=1$,
$\Omega_k|_{\mathbb R\setminus f(\gamma_k,t)}=0$, $k=1,2$. The solution of the problem is known,
see, e.g., [5,~p.334]. Namely, $$m_k(t)=\frac{\alpha_k(t)}{\pi}$$ where $\alpha_k(t)$ is the angle
under which the segment $I_k=I_k(t)$ is observed from the point $w=i$, $k=1,2$. It remains to find
asymptotic expansions for $\alpha_k(t)$.

Since $$f_1(0,t)=\root3\of
t+6\root3\of{t^2}+O(t),\;\;\;f_2(0,t)=-3\root3\of{t^2}+O(t),\;\;\;t\to+0,$$ after elementary
geometrical considerations we have $$\alpha_1(t)=\arctan f_1(0,t)-\arctan\root3\of
t=6\root3\of{t^2}+O(t),\;\;\;t\to+0,$$ $$\alpha_2(t)=\arctan\root3\of t-\arctan f_2(0,t)=\root3\of
t+3\root3\of{t^2}+O(t),\;\;\;t\to+0.$$ This implies that
$$\frac{m_1(t)}{m_2^2(t)}=\pi\frac{6\root3\of{t^2}+O(t)}{(\root3\of t+3\root3\of{t^2}+O(t))^2}=
6\pi(1+O(\root3\of t)),\;\;\;t\to+0,$$ which leads to (\ref{har}) and completes the
proof.
\endproof

\begin{remark}
The relation similar to (\ref{har}) follows from \cite{ProkhVas} for the two sides of the circular
slit $\gamma(t)$ in $\mathbb H$ such that $\gamma(t)$ is tangential to $\mathbb R$ at $z=0$.
\end{remark}

\section{Representation of holomorphic solutions}

Holomorphic solutions to (\ref{Le3}) or, equivalently, to (\ref{Le2}) appear in a neighborhood of
every non-singular point $(z_0,0)$. We will be interested in real solutions corresponding to
$z_0\in\mathbb R$.

Put $z_0=\epsilon>0$ and let

\begin{equation}
f(\epsilon,t)=\epsilon+\sum_{n=1}^{\infty}a_n(\epsilon)t^{n/3} \label{hol}
\end{equation}
be a solution of equation (\ref{Le3}) holomorphic with respect to $\tau=\root3\of t$. Change
$\root3\of t$ by $\tau$ and substitute (\ref{hol}) in (\ref{Le2}) to get that

\begin{equation}
\sum_{n=1}^{\infty}na_n(\epsilon)\tau^{n-1}\left[\epsilon-\tau+\sum_{n=1}^{\infty}
a_n(\epsilon)\tau^n\right]=6\tau^2. \label{coe}
\end{equation}

Equate coefficients at the same powers in both sides of (\ref{coe}) and obtain equations

\begin{equation}
a_1(\epsilon)=0,\;\;\;a_2(\epsilon)=0,\;\;\;a_k(\epsilon)=\frac{6}{k\epsilon^{k-2}}, \;\;\;k=3,4,5,
\label{de1}
\end{equation}
and
\begin{equation}
a_n(\epsilon)=\frac{1}{n\epsilon}\left[(n-1)a_{n-1}(\epsilon)-
\sum_{k=3}^{n-3}(n-k)a_{n-k}(\epsilon)a_k(\epsilon)\right], \;\;\;n\geq6. \label{de2}
\end{equation}

The series in (\ref{hol}) converges for $|\tau|=|\root3\of t|<R(\epsilon)$.

\begin{theorem}
The series in (\ref{hol}) converges for
\begin{equation}
|t|<\epsilon^3+o(\epsilon^3),\;\;\;\epsilon\to+0. \label{rad}
\end{equation}
\end{theorem}

\proof Estimate the convergence radius $R(\epsilon)$ following the Cauchy majorant method, see,
e.g., [4,~Chap.1, \S\S2-3], [16,~Chap.3, \S1]. The Cauchy theorem states: if the right-hand side in
(\ref{Le2}) is holomorphic on a product of the closed disks $|g-\epsilon|\leq\rho_1$ and
$|\tau|\leq r_1$ and is bounded there by $M$, then the series
$\sum_{n=1}^{\infty}a_n(\epsilon)\tau^n$ converges in the disk
$$|\tau|<R(\epsilon)=r_1\left(1-\exp\left\{-\frac{\rho_1}{2Mr_1}\right\}\right).$$

In the case of equation (\ref{Le2}) we have
$$\rho_1+r_1<\epsilon,\;\;\text{and}\;\;M=\frac{6r_1^2}{\epsilon-(\rho_1+r_1)}.$$ This implies
that for $\rho_1+r_1=\epsilon-\delta,$ $\delta>0$,
$$R(\epsilon)=r_1\left(1-\exp\left\{-\frac{\epsilon-
\delta-r_1}{12r_1^2}\delta\right\}\right).$$

So $R(\epsilon)$ depends on $\delta$ and $r_1$. Maximum of $R$
with respect to $\delta$ is obtained for
$\delta=(\epsilon-r_1)/2$. Hence, this maximum is equal to
\begin{equation}
R_1(\epsilon)=r_1\left(1-\exp\left\{-\frac{(\epsilon-r_1)^2}{48r_1^3}\right\}\right),
\label{eR1}
\end{equation}
where $R_1(\epsilon)$ depends on $r_1$. Let us find a maximum of $R_1$ with respect to
$r_1\in(0,\epsilon)$. Notice that $R_1$ vanishes for $r_1=0$ and $r_1=\epsilon$. Therefore the
maximum of $R_1$ is attained for a certain root $r_1=r_1(\epsilon)\in(0,\epsilon)$ of the
derivative of $R_1$ with respect to $r_1$. To simplify the calculations we put
$r_1(\epsilon)=\epsilon c(\epsilon)$, $0<c(\epsilon)<1$. Now the derivative of $R_1$ vanishes for
$c=c(\epsilon)$ satisfying
\begin{equation}
1-\exp\left\{-\frac{(1-c)^2}{48\epsilon c^3}\right\}\left(1+\frac{(1-c)(3-c)}{48\epsilon
c^3}\right)=0. \label{der}
\end{equation}

Choose a sequence $\{\epsilon_n\}_{n=1}^{\infty}$ of positive
numbers, $\lim_{n\to\infty}\epsilon_n=0$, such that
$c(\epsilon_n)$ converge to $c_0$ as $n\to\infty$. Suppose that
$c_0<1$. Then
$$\exp\left\{-\frac{(1-c(\epsilon_n))^2}{48\epsilon_nc^3(\epsilon_n)}\right\}
\left(1+\frac{(1-c(\epsilon_n))(3-c(\epsilon_n))}{48\epsilon_nc^3(\epsilon_n)}\right)<1$$
for $n$ large enough. Therefore $c(\epsilon_n)$ is not a root of
equation (\ref{der}) for $\epsilon=\epsilon_n$ and $n$ large
enough. This contradiction  claims that $c_0=1$ for every sequence
$\{\epsilon_n>0\}_{n=1}^{\infty}$ tending to 0 with
$\lim_{n\to\infty}c(\epsilon_n)=c_0$. So we proved that
$c(\epsilon)\to1$ as $\epsilon\to+0$.

Consequently, the maximum of $R_1$ with respect to $r_1$ is attained for
$r_1(\epsilon)=\epsilon c(\epsilon)=\epsilon(1+o(1))$ as $\epsilon\to+0$. Let
$R_2=R_2(\epsilon)$ denote the maximum of $R_1$ with respect to $r_1$. It follows from
(\ref{eR1}) that
\begin{equation}
R_2(\epsilon)=r_1(\epsilon)
\left(1-\exp\left\{-\frac{(\epsilon-r_1(\epsilon))^2}{48r_1^3(\epsilon)}\right\}\right)=
\epsilon c(\epsilon)\left(1-\exp\left\{-\frac{(1-c(\epsilon))^2}{48\epsilon
c^3(\epsilon)}\right\}\right). \label{max}
\end{equation}

Examine how fast does $c(\epsilon)$ tends to 1 as $\epsilon\to+0$. Choose a sequence
$\{\epsilon_n>0\}_{n=1}^{\infty}$, $\lim_{n\to\infty}\epsilon_n=0$, such that the sequence
$(1-c(\epsilon_n))^2/\epsilon_n$ converges to a non-negative number or to $\infty$. Denote
$$l:=\lim_{n\to\infty}\frac{(1-c(\epsilon_n))^2}{\epsilon_n},\;\;\;0\leq l\leq\infty.$$

If $0<l<\infty$, then $(1-c(\epsilon_n))/\epsilon_n$ tends to $\infty$, and equation
(\ref{der}) with $\epsilon=\epsilon_n$ has no roots for $n$ large enough.

If $l=0$, then, according to (\ref{der}), $\lim_{n\to\infty}(1-c(\epsilon_n))/\epsilon_n=0$,
and $$\exp\left\{-\frac{(1-c(\epsilon_n))^2}{48\epsilon
c^3(\epsilon)}\right\}\left(1+\frac{(1-c(\epsilon_n))(3-c(\epsilon_n))}{48\epsilon_n
c^3(\epsilon_n)}\right)=$$ $$\left(1-\frac{(1-c(\epsilon_n))^2}{48\epsilon_n
c^3(\epsilon_n)}+o\left(\frac{(1-c(\epsilon_n))^2}{\epsilon_n}\right)\right)
\left(1+\frac{(1-c(\epsilon_n))(3-c(\epsilon_n))}{48\epsilon_n c^3(\epsilon_n)}+\right)=$$
$$1+\frac{1-c(\epsilon_n)}{24\epsilon_n}+ o\left(\frac{1-c(\epsilon_n)}{\epsilon_n}\right),
\;\;\;n\to\infty.$$ This implies again that equation (\ref{der}) with $\epsilon=\epsilon_n$ has no
roots for $n$ large enough.

Thus the only possible case is $l=\infty$ for all sequences $\{\epsilon_n>0\}_{n=1}^{\infty}$
converging to 0. It follows from (\ref{max}) that
\begin{equation}
R_2(\epsilon)=\max_{0<r_1(\epsilon)<\epsilon}R_1(\epsilon)=\epsilon+o(\epsilon),
\;\;\;\epsilon\to0. \label{R2}
\end{equation}
In other words, the series in (\ref{hol}) converges for $|t|<(\epsilon+o(\epsilon))^3$,
$\epsilon\to0$, which implies the statement of Theorem 3 and completes the proof.
\endproof

\begin{remark}
Evidently, a similar conclusion with the same formulas (\ref{de1}) and (\ref{de2}) is true for
$\epsilon<0$.
\end{remark}

\end{document}